\input amstex
\documentstyle{amsppt}

\centerline{\bf ON $\bold{5}$-CYCLES AND STRONG $\bold{5}$-SUBTOURNAMENTS }

\smallskip

\centerline{\bf IN A TOURNAMENT OF ODD ORDER $\bold{N}$}

\bigskip

\centerline{\bf S.V. Savchenko}

\bigskip

\centerline{ L.D. Landau Institute for Theoretical Physics, Russian Academy of Sciences}
\smallskip
\centerline{ Kosygin str. 2, Moscow 119334, Russia
\footnote[]{E-mail: savch$\@$itp.ac.ru}}

\bigskip

\centerline{\sl Dedicated to Richard Brualdi on the occasion
of his 80th birthday}

\bigskip

{\bf Abstract.} Let $T$ be a tournament of odd order $n\ge 5,$ $c_{m}(T)$
be the number of its $m$-cycles,
and $s_{m}(T)$ be the number of its strongly connected $m$-subtournaments.
Due to work of L.W. Beineke and F. Harary,
it is well known that
$s_{m}(T)\le s_{m}(RLT_{n}),$
where $RLT_{n}$ is the regular locally transitive tournament of order $n.$
For $m=3$ and $m=4,$ $c_{m}(T)$ equals $s_{m}(T),$ but it is not so
for $m\ge 5.$ As J.W. Moon pointed out in his note in 1966,
the problem of determining the maximum of $c_{m}(T)$ seems very difficult
in general (i.e. for $m\ge 5$). In the present paper,
based on the Komarov-Mackey formula for $c_{5}(T)$ obtained recently,
we prove that
$c_{5}(T)\le (n+1)n(n-1)(n-2)(n-3)/160$ with equality holding iff $T$ is doubly regular.
A formula for $s_{5}(T)$ is also deduced. With the use of it,
we show that
$s_{5}(T)\le (n+1)n(n-1)(n-3)(11n-47)/1920$ with equality holding iff
$T=RLT_{n}$ or $n=7$ and $T$ is regular or $n=5$ and $T$ is strong.
It is also proved that for a regular tournament $T$
of (odd) order $n\ge 9,$ a lower bound
$(n+1)n(n-1)(n-3)(17n-59)/3840\le s_{5}(T)$
holds with equality iff $T$ is doubly regular.
These results are compared with the ones
recently obtained by the author for $c_{5}(T).$

\bigskip


\bigskip


\bigskip

\centerline{\bf \S 1. Introduction}

\medskip

A {\sl tournament} $T$ of order $n$ (or, merely, $n$-tournament $T$)
is an orientation of the complete
graph $K_{n}.$ So, there exists exactly one arc between any two vertices
of $T.$ If a pair $(i,j)$ is an arc in $T,$ we say that the vertex $i$
{\sl dominates} the vertex $j$ and write $i\to j.$
For two vertex-disjoint
sets (tournaments) $I$ and $J,$ the notation $I\Rightarrow J$ means
that each vertex of $I$ dominates every vertex of $J.$
We say that vertices of $T$ are replaced by some tournaments $S^{(1)},...,S^{(n)}$
if the latter are taken instead of the former and the binary
relation $\to$ between the vertices is replaced by the binary relation
$\Rightarrow$ between the tournaments. For a tournament $S$,
we write $T\cdot S$ if each
vertex of $T$ is replaced by a copy of $S.$ The tournament obtained in
result is called the {\sl composition} of $T$ and $S.$

Let $N^{+}(i)$ be the {\sl out-set} of $i$ (i.e. the set of vertices
dominated by $i$) and
$N^{-}(i)$ be the {\sl in-set} of $i$
(i.e. the set of vertices dominating $i$).
The subtournaments induced by these vertex-sets will be denoted by the
same symbols.
The quantities $\delta^{+}_{i}=|N^{+}(i)|$ and $\delta^{-}_{i}=|N^{-}(i)|$
are called the {\sl out-degree} and {\sl in-degree} of $i,$ respectively.
They are related as follows:
$\delta^{+}_{i}=n-1-\delta^{-}_{i}.$
However, we will consider both
the out-degrees and the in-degrees of the vertices,
which are also called
the {\sl scores} and {\sl co-scores} of the vertices.
The {\sl score-list} ({\sl co-score-list}) of
a tournament is the list of the out-degrees (in-degrees) of
the vertices, usually arranged in non-decreasing order.
A score-list is {\sl balanced} if it coincides with the
co-score-list. We also say that a non-decreasing
sequence of non-negative integers
$\delta_{1}\le ... \le \delta_{k}\le ... \le \delta_{n}$
is a (co)-score-list if there is a tournament of order $n$
whose (co)-scores are $\delta_{1},...,\delta_{k},...,\delta_{n}.$
According to the Landau criterion (see [21], [26], and [9]), this holds
if and only if $\delta_{1}+...+\delta_{k}\ge \binom{k}{2}$
for each $k=1,...,n-1$ and $\delta_{1}+...+\delta_{n}=\binom{n}{2}.$

A vertex in $T$ is a {\sl source} if its out-degree is $n-1.$
In turn, a vertex in $T$ is a {\sl sink} if its in-degree is $n-1.$
Any $n$-tournament admits at most one source
and at most one sink. Moreover, there exists precisely one $n$-tournament
each of whose subtournaments has exactly one sink and exactly one source,
namely, the {\sl transitive} tournament $TT_{n}$ of order $n.$
It is called so because if $i\to k\to j,$ then $i\to j$ and this (transitive)
condition uniquely determines $TT_{n}$: if $1,...,n$ is its
(unique) Hamiltonian path, then $j\to i$ if and only if $i>j.$

A tournament is {\sl strongly connected} (or, merely, {\sl strong})
if for any two distinct vertices $i$ and $j,$ there is a path from
$i$ to $j.$ We assume that
the tournament $\circ$ of order $1$ (consisting of a unique vertex)
is strong. Denote by $\Cal{S}\Cal{T}_{n}$ the set of all strong tournaments
of order $n.$
If $T$ is not contained in $\Cal{S}\Cal{T}_{n}$, then
it is obtained from some $TT_{k},$
where $2\le k\le n,$ by replacing its vertices $1,...,k$ with strong tournaments
$T^{(1)},...,T^{(k)},$ respectively.
In this case, we write
$T=T^{(1)}\Rightarrow...\Rightarrow T^{(k)}.$

For a family $\Cal{F}_{m}$ of digraphs of order $m,$
let $n_{\Cal{F}_{m}}(T)$ be the number of
copies of elements of $\Cal{F}_{m}$ in $T$ (i.e. the number of subdigraphs of $T$ that are
isomorphic to elements of $\Cal{F}_{m}$).
If $\Cal{F}_{m}$ consists of the unique element $D_{m},$ then we simply write
$n_{D_{m}}(T).$
A problem of determining the maximum (minimum) of
$n_{\Cal{F}_{m}}(T)$ or $n_{D_{m}}(T)$
in the class ${\Cal T}_{n}$ of all tournaments of order $n$ or
some subclass of ${\Cal T}_{n}$ can be posed. In the present paper,
we consider the cases of $\Cal{F}_{m}=\Cal{S}\Cal{T}_{m}$ and $D_{m}=\vec{C}_{m},$
where $\vec{C}_{m}$ is the {\sl directed cycle of length $m$}.
For simplicity, denote $n_{\Cal{S}\Cal{T}_{m}}(T)$ and $n_{\vec{C}_{m}}(T)$
by $s_{m}(T)$ and $c_{m}(T),$ respectively.

We first consider the quantity $s_{m}(T).$
As for $m\ge 3,$ we have $\Cal{S}\Cal{T}_{m}\cap
\Cal{T}_{m-1}\Rightarrow\circ=\emptyset,$
the inequality
$$s_{m}(T)
\le \binom{n}{m}-n_{\Cal{T}_{m-1}\Rightarrow\circ}(T)\eqno(1)$$
holds with equality iff every subtournament of order $m$ in $T$ is either
strong or contains a sink.
In particular, this condition is satisfied if $T$
is {\sl locally$^{+}$ transitive}, i.e. the
out-set of each vertex induces a transitive subtournament.
Note that in $(1),$ the term $n_{\Cal{T}_{m-1}\Rightarrow\circ}(T)$ can be replaced by
$n_{\circ\Rightarrow\Cal{T}_{m-1}}(T).$
The new upper bound will be attained when $T$
is {\sl locally$^{-}$ transitive}, i.e. the
in-set of each vertex induces a transitive subtournament.
As
$$\sum\limits_{i=1}^{n}n_{TT_{3}}\Bigl(N^{+}(i)\Bigr)=
n_{\circ\Rightarrow TT_{3}}(T)=n_{TT_{4}}(T)=
n_{TT_{3}\Rightarrow \circ}(T)=
\sum\limits_{i=1}^{n}n_{TT_{3}}\Bigl(N^{-}(i)\Bigr),$$
$T$ with balanced (co)-score-list is locally$^{+}$ transitive
if and only if it is locally$^{-}$ transitive.
A tournament having both of these properties is {\sl locally transitive}.
According to Theorem 4 [1], for a given balanced
(co)-score-list $\delta_{1}\le ... \le \delta_{n},$ there exists a locally transitive tournament with this
(co)-score-list if and only if any integer between $\delta_{1}$ and $\delta_{n}$
is contained in $\{\delta_{1},...,\delta_{n}\}.$

Note that
$$n_{\Cal{T}_{p}\Rightarrow\circ}(T)=
\sum\limits_{i=1}^{n}n_{\Cal{T}_{p}}
\Bigl(N^{-}(i)\Bigr)=
\sum\limits_{i=1}^{n}\binom{\delta_{i}^{-}}{p}.\eqno(2)$$
A well-known combinatorial result (presented, for instance, in [4])
states that for a given integer $p\ge 2,$
the sum is a minimum if
$\delta_{i}^{-}$ are as nearly equal as possible
(the same also holds if one replaces $\delta_{i}^{-}$ by $\delta_{i}^{+}$
in the sum). Moreover,
by our Lemma 1, for $n\ge 2p-1,$ "if" can be replaced by "if and only if"
in this statement. The condition that $\delta_{i}^{-}$ are as nearly equal
as possible means that if $n$ is odd, each $\delta_{i}^{-}$ equals $\frac{n-1}{2}$ and
hence, $T$ is {\sl regular} (for this case,
each $\delta_{i}^{+}$ also equals $\frac{n-1}{2};$ so, this quantity
can be called the {\sl semi-degree} and denoted by $\delta$);
if $n$ is even, half the in-degrees
are $\frac{n}{2}$ and the others are $\frac{n}{2}-1,$ i.e. $T$
is {\sl near regular} (as one can see,
it is a one-vertex-deleted subtournament of some regular tournament
of order $n+1$).

It is not difficult to check that for each odd $n,$
there exists exactly one $n$-tournament that is both regular and
locally transitive (see [2]).
It can be defined on
the ring ${\Bbb Z}_{n}$ of residues modulo $n$ by the rule $i\to j$ iff
the difference $j-i$ (as a residue) is contained in the subset
$\{1,...,\frac{n-1}{2}\}$ of ${\Bbb Z}_{n}.$
We call it
the {\sl regular locally transitive} tournament and denote it by $RLT_{n}.$
The above arguments imply that for odd $n,$
the maximum of $s_{m}(T)$ in the class ${\Cal T}_{n}$
is attained at $RLT_{n}$ (see [4] and [10]).
However, to determine all maximizers of $s_{m}(T)$ in the class
${\Cal T}_{n},$ one needs to get an exact expression for $s_{m}(T).$
It is possible to do when $m=3,4,$ and $5$ because for these values of $m,$
any non-strong tournament of order $m$ admits either a sink or a source. So,
in these cases, we have
$$s_{m}(T)=\binom{n}{m}-n_{\Cal{T}_{m-1}\Rightarrow\circ}(T)-n_{\circ\Rightarrow \Cal{T}_{m-1}}(T)
+n_{\circ\Rightarrow \Cal{T}_{m-2}\Rightarrow\circ}(T)=$$
$$\binom{n}{m}-\sum\limits_{i=1}^{n}\binom{\delta_{i}^{-}}{m-1}
-\sum\limits_{i=1}^{n}\binom{\delta_{i}^{+}}{m-1}+
\sum\limits_{i=1}^{n}n_{\Cal{T}_{m-2}\Rightarrow\circ}
\Bigl(N^{+}(i)\Bigr).\eqno(3)$$

Note that $\circ\Rightarrow\circ=TT_{2}.$
As each $2$-subtournament of any tournament is $TT_{2}$,
for $m=3,$ the last two terms in $(3)$ are cancelled in pairs.
So, the maximum of $s_{3}(T)$ in the class ${\Cal T}_{n}$
is attained only at the class $\Cal{R}_{n}$
of regular tournaments of order $n$ or the class $\Cal{NR}_{n}$
of near regular tournaments of order $n$ if $n$ is odd or even, respectively [15].
In turn, as $\Cal{T}_{2}\Rightarrow\circ=\{TT_{2}\Rightarrow\circ\}=\{TT_{3}\}$
and a tournament admits only transitive $3$-subtournaments
if and only if it is transitive itself,
for each odd $n\ge 5,$ the maximum of $s_{4}(T)$ in the class ${\Cal T}_{n}$
is attained if and only if $T$ is regular and locally transitive, i.e. $T=RLT_{n}.$
Finally, one can show that
a tournament of order $p\ge 4$ admits only $4$-subtournaments
belonging to $\Cal{T}_{3}\Rightarrow\circ$
if and only if it is contained in $\Cal{T}_{3}\Rightarrow TT_{p-3}.$
Based on this fact and $(3),$ we prove in Section 2
that for each odd $n\ge 9,$ the maximum of $s_{5}(T)$ in the class
$\Cal{T}_{n}$ is attained only at $RLT_{n},$ while for $n=7,$ it is
achieved if and only if $T$ is regular.

For $m=4$ and $m=5,$ the quantity $s_{m}(T)$ can take different values on elements of $\Cal{R}_{n}$
or $\Cal{NR}_{n}.$
A lower bound on $s_{m}(T)$ in the class of tournaments of order $n$ with given
(co)-score-list can be obtained from $(3)$ if
for each $i=1,...,n,$ one replaces $n_{\Cal{T}_{m-2}\Rightarrow\circ}\Bigl(N^{+}(i)\Bigr)$
by its minimum possible value. (For $m=4,$ this was first done in [18].)
So, equality $(2)$ and the above lemma on the combinatorial sum imply that
the lower bound on $s_{m}(T)$ is attained if for each $i=1,...,n,$
the out-set $N^{+}(i)$ induces a regular
or near regular tournament
when $\bigl|N^{+}(i)\bigr|$ is odd or even, respectively.
A tournament having this property is called {\sl locally$^{+}$ regular}.

In the last sum in $(3),$ the term
$n_{\Cal{T}_{m-2}\Rightarrow\circ}\Bigl(N^{+}(i)\Bigr)$ can be replaced by
the term $n_{\circ\Rightarrow\Cal{T}_{m-2}}\Bigl(N^{-}(i)\Bigr).$
So, a new lower bound on $s_{m}(T)$ can be obtained if
for each $i=1,...,n,$ one replaces $n_{\circ\Rightarrow\Cal{T}_{m-2}}\Bigl(N^{-}(i)\Bigr)$
by its minimum possible value. It is attained if for each $i=1,...,n,$
the in-set $N^{-}(i)$ induces a regular
or near regular tournament
when $\bigl|N^{-}(i)\bigr|$ is odd or even, respectively.
Such a tournament is called {\sl locally$^{-}$ regular}.
One can show that a tournament with balanced (co)-score-list is locally$^{+}$ regular
if and only if it is locally$^{-}$ regular.
In this case, it is {\sl locally regular}.

A regular locally regular tournament of order $n$
is called {\sl doubly-regular} or {\sl nearly-doubly-regular}
when $\frac{n-1}{2}$ is odd and hence, $n\equiv3\mod 4$ or
$\frac{n-1}{2}$ is even and hence, $n\equiv1\mod 4$, respectively.
Denote by $\Cal{DR}_{n}$
and $\Cal{NDR}_{n}$ the corresponding subfamilies
of $\Cal{R}_{n}.$ If the value of $s_{m}(T)$ does not depend on a particular
choice of $T\in \Cal{DR}_{n}$ or $T\in \Cal{NDR}_{n}$, then we will simply write
$s_{m}(\Cal{DR}_{n})$ or $s_{m}(\Cal{NDR}_{n}).$ The same rule also concerns
the other quantities introduced below.

Under this convention,
the lower bound on $s_{m}(T)$ in the class $\Cal{R}_{n}$ suggested above for
the case of $m=4$ and $m=5$  can be written as $s_{m}\bigl(\Cal{(N)DR}_{n}\bigr) \le s_{m}(T),$
where the choice of $\Cal{NDR}_{n}$ or $\Cal{DR}_{n}$ depends on
the residue ($1$ or $3$) of $n$  modulo $4.$
Our Lemma 1 implies that for $m=4$ and $n\ge 5$ or $m=5$ and $n\ge 11,$
the equality holds
if and only if $T\in \Cal{DR}_{n}$ (when $n\equiv 3\mod 4$) or $T\in \Cal{NDR}_{n}$ (when $n\equiv 1\mod 4$).
However, it is not so for $m=5$ and $n=9.$
For this case, the bound is attained if and only if the out-set of each vertex of $T$ includes no sink.
In Section 2, we show that besides two elements of $\Cal{NDR}_{9}$,
there exist three more such regular tournaments of order $9.$

The problem of non-emptiness of $\Cal{DR}_{n}$ for each
$n\equiv 3\mod 4$ is open up to now, while one can present infinitely
many $DR$'s, there are methods for constructing $DR$'s
from those of smaller orders and according to a common opinion,
at least one element of $\Cal{DR}_{n}$ exists for any possible order $n$ (see [25]).
The same can be also said about
$\Cal{NDR}_{n}$, where $n\equiv1\mod 4,$ (see [22], [23], and [33]),
while much less is known about $NDR$'s than about $DR$'s.
We say that for odd $n,$ the
$\exists$-property holds
if the set $\Cal{RLR}_{n}$ of regular locally regular tournaments of order $n$ is not empty,
i.e. either $\Cal{DR}_{n}\neq \emptyset$ (when $n\equiv3\mod 4$)
or $\Cal{NDR}_{n}\neq \emptyset$ (when $n\equiv1\mod 4$).

Note that
if $m=3$ or $4,$ then there exists exactly one strong tournament
of order $m.$
For $m=3,$ it is called the {\sl cyclic triple} and is denoted by
$\Delta.$ In turn, the unique strong tournament
$ST_{4}$ of order $4$ can be obtained from
$\Delta$ by replacing one of its vertices with $TT_{2}.$
It is not difficult to check that either of the strong tournaments
contains precisely one Hamiltonian cycle.
So, for $m=3$ and $m=4,$  the equality $c_{m}(T)=s_{m}(T)$ holds
and hence, we can apply the above-mentioned (classical)
results on $s_{m}(T)$ to $c_{m}(T).$

As it was first pointed out in [20] (see also [21]), for the case of $m=5,$
the problem of determining the maximum of $c_{m}(T)$ in the class $\Cal{T}_{n}$
is much more difficult.
(Note that $|\Cal{S}\Cal{T}_{5}|=6$ and the number of $5$-cycles
varies from $1$ to $3$ on $\Cal{S}\Cal{T}_{5}$.)
One can show that in contrast with $s_{5}(T)$ and $c_{4}(T),$
the collection of the out-sets of vertices in $T$ does not determine $c_{5}(T).$
However, there also exists a formula for $c_{5}(T)$ expressed
in terms of the co-scores of vertices in
the out-sets of vertices in $T.$
It has been obtained only recently
by Komarov and Mackey (see [16]). With the use of it,
they present an upper bound on
$c_{5}(T)$ in the class ${\Cal T}_{n}$ (which, however, is not sharp for
any $n\ge 5$) and prove that
the maximum of $c_{5}(T)$ in the class ${\Cal T}_{n}$ is asymptotically the
same as the expected number $E_{n}c_{5}$ of $5$-cycles in a random tournament
of order $n.$ In Section 3, based on the Komarov-Mackey formula,
we show that for each $T\in {\Cal T}_{n},$
where $n$ is odd and $n\ge 5,$
the inequality $c_{5}(T)\le E_{n+1}c_{5}$ holds with equality
iff $T\in \Cal{DR}_{n}.$ (In particular, this means that the upper bound is sharp for infinitely
many $n$.) So, if the (very plausible) conjecture on the existence of
at least one element in $\Cal{DR}_{n}$ for each
possible $n$ is true, then the problem of determining the maximum of $c_{5}(T)$
in the class ${\Cal T}_{n}$ is settled for $n\equiv 3\mod 4.$
The other cases are still open.
At the end of the section, the reader will find a conjecture on possible
maximizers of $c_{5}(T)$ in ${\Cal T}_{n}$ for $n\equiv 1\mod 4.$
Finally, in Section 4, we give concluding comments on $s_{m}(T)$ and $c_{m}(T)$
in the case of arbitrary $m.$

\bigskip
\centerline{\bf \S 2. All maximizers of $\bold{s_{5}(T)}$
in the class $\bold{{\Cal T}_{n}},$ where $\bold{n}$ is odd,}
\smallskip
\centerline{\bf and a lower bound on $\bold{s_{5}(T)}$ in the class $\bold{{\Cal R}_{n}}$}
\medskip

We have already mentioned and used the well-known fact that the minimum of the sum
$\sum\limits_{i=1}^{n}\binom{\delta_{i}^{-}}{p}$
is achieved when $\delta_{i}^{-}$ are as nearly equal as possible.
For our purposes, we need to show that this condition is also necessary
if $p\ge 2$ and $n$ is large enough.

\smallskip

{\bf Lemma 1.} {\sl Let $p\ge 2.$
Then for $n\ge 2p-1,$ the sum
$\sum\limits_{i=1}^{n}\binom{\delta_{i}^{-}}{p}$
attains its minimum in the class $\Cal{T}_{n}$
if and only if the corresponding tournament is
regular (when $n$ is odd) or near regular (when $n$ is even).}

\smallskip
{\bf Proof.} We first consider the minimization problem
in the class $\Omega_{n}$
of ordered sequences $\omega_{n}$
of non-negative integers $\delta_{1}^{-}\le ...\le\delta_{n}^{-}$ whose sum is equal to
$\frac{n(n-1)}{2}$ for the case of $n\ge 2p.$
If $\delta_{1}^{-}\ge p-1$ and
$\delta_{n}^{-}> \delta_{1}^{-}+1,$ then
$$\binom{\delta_{1}^{-}+1}{p}-\binom{\delta_{1}^{-}}{p}=
\frac{\delta_{1}^{-}...(\delta_{1}^{-}-p+2)}{(p-1)!}<$$
$$\frac{(\delta_{n}^{-}-1)...(\delta_{n}^{-}-p+1)}{(p-1)!}=
\binom{\delta_{n}^{-}}{p}-\binom{\delta_{n}^{-}-1}{p}.$$
In turn, if $\delta_{1}^{-} < p-1,$ then this inequality also holds
because the LHS equals $0,$
while the RHS is positive as $\delta_{n}^{-}\ge\lceil\frac{n-1}{2}\rceil \ge p.$
So, in both cases, we have
$$\binom{\delta_{1}^{-}}{p}+
\binom{\delta_{n}^{-}}{p}> \binom{\delta_{1}^{-}+1}{p}+\binom{\delta_{n}^{-}-1}{p}$$
and hence,
the sum $\sum\limits_{i=1}^{n}\binom{\delta_{i}^{-}}{p}$ will strictly decrease if we replace
$\delta_{n}^{-}$ and $\delta_{1}^{-}$ with $\delta_{n}^{-}-1$ and
$\delta_{1}^{-}+1,$ respectively.
(Obviously, after ordering, the new sequence will belong to $\Omega_{n}$.)

We see that the minimum is achieved if $\delta_{n}^{-}-\delta_{1}^{-}\le 1.$
Let us show that this condition uniquely determines $\omega_{n}$ in the class
$\Omega_{n}.$ If $\delta_{n}^{-}-\delta_{1}^{-}=0,$ then
$\delta_{1}^{-}=...=\delta_{n}^{-}=\frac{n-1}{2},$ where $n$ must be odd.
For any other $\omega_{n},$ we have $\delta_{1}^{-}<\frac{n-1}{2}<\delta_{n}^{-}$
and hence, if $n$ is odd, then $\delta_{n}^{-}-\delta_{1}^{-}\ge 2.$
So, if $\delta_{n}^{-}-\delta_{1}^{-}=1,$ then $n$ must be even,
$\delta_{1}^{-}=\frac{n}{2}-1$ and $\delta_{n}^{-}=\frac{n}{2}.$ Let $p$ be the integer such that
$\frac{n}{2}-1=\delta_{p}^{-} < \delta_{p+1}^{-}=\frac{n}{2}.$
Then the condition $\delta_{1}^{-}+...+\delta_{n}^{-}=p\bigl(\frac{n}{2}-1\bigr)
+(n-p)\frac{n}{2}=\frac{n(n-1)}{2}$
implies that $p=\frac{n}{2}.$
In both cases, denote the sequence obtained in result by $\hat{\omega}_{n}.$

The arguments presented above imply that for $n\ge 2p,$
the minimum of the sum $\sum\limits_{i=1}^{n}\binom{\delta_{i}^{-}}{p},$ where $p\ge 2,$
in $\Omega_{n}$ is attained only at $\hat{\omega}_{n}.$
In the case where $n=2p-1,$
the sequence $\hat{\omega}_{n}$ is a unique element of $\Omega_{n}$ for which
the sum is equal to zero.
(For any other sequence of $\Omega_{n},$ where $n=2p-1,$ we have
$\delta_{n}^{-}\ge p$ and hence, the sum is greater than or equal to $1.$)
So, we have determined all minimizers of the sum in $\Omega_{n}$ for any
$n\ge 2p-1.$ The same result also holds for the minimum in the class $\Cal{T}_{n}$
because $\hat{\omega}_{n}$ is the co-score-list of a regular $n$-tournament
or a near regular $n$-tournament when $n$ is odd or even, respectively.
The lemma is proved.

\smallskip

Denote by $\Delta\cdot TT_{2}$ the (near regular) tournament of order $6$
obtained from $\Delta$ by replacing each vertex with $TT_{2}.$
It is well known that
there exist exactly three regular tournaments of order $7,$
namely, $RLT_{7},$ the unique element $DR_{7}$ of $\Cal{DR}_{7},$
and the Kotzig tournament $Kz_{7}$
which is uniquely determined by the condition that $\Delta\cdot TT_{2}$
is its one-vertex-deleted subtournament.
As we have already seen, for each odd $n,$
the maximum of $s_{5}(T)$ in the class $\Cal{T}_{n}$
is equal to
$$s_{5}(RLT_{n})=\binom{n}{5}-n\binom{\frac{n-1}{2}}{4}=
\frac{(n+1)n(n-1)(n-3)(11n-47)}{1920}.\eqno(4)$$
The description of $\Cal{R}_{7}$ and
Lemma 1 taken together allow us to determine
all maximizers of $s_{5}(T)$ in the class $\Cal{T}_{n}.$

\smallskip

{\bf Proposition 1.} {\sl For each odd $n\ge 9,$  the maximum of $s_{5}(T)$
in the class $\Cal{T}_{n}$ is attained only at $RLT_{n}.$
For $n=7,$ it is also achieved for $DR_{7}$ and $Kz_{7}.$}

\smallskip
{\bf Proof.} Recall that $\Cal{T}_{4}$ consists of $ST_{4},$ $\circ\Rightarrow\Delta,$
$\Delta\Rightarrow\circ,$ and $TT_{4}.$  As
$\{\Delta\Rightarrow\circ\} \cup \{TT_{4}\}=\Cal{T}_{3}\Rightarrow\circ,$
for $m=5,$ we can rewrite $(3)$ as
$$s_{5}(T)=
\binom{n}{5}-\sum\limits_{i=1}^{n}\binom{\delta^{-}_{i}}{4}
-\sum\limits_{i=1}^{n}n_{ST_{4}}\Bigl(N^{+}(i)\Bigr)-
\sum\limits_{i=1}^{n}n_{\circ\Rightarrow\Delta}\Bigl(N^{+}(i)\Bigr).
\eqno(5)$$
For the case of $RLT_{n},$ each of the three sums in $(5)$ attains its minimum.
The first one equals $n\binom{\frac{n-1}{2}}{4}$ and the others
equal $0$ (see $(4)$).
The same must also hold for
any maximizer $T$ of order $n\ge 7=2\cdot 4-1.$
Equality $(5)$ and Lemma 1 imply that $T$ is a regular tournament
with semi-degree $\delta\ge 3$ the out-set of each of whose vertices
contains no $ST_{4}$ or $\circ\Rightarrow\Delta$ as a subtournament.
Obviously, this condition always holds if $\delta=3.$ So, we can
assume that $\delta \ge 4.$

Any tournament of order $\delta$ is obtained from some transitive tournament
by replacing its vertices with
strong tournaments. The latter are strong components of the tournament.
If the order of at least one of them is not less than $4,$ then
by the Moon vertex-pancyclic theorem, one can always find a copy of $ST_{4}$ in it.
So, if the tournament contains no $ST_{4},$ then each of its strong components
is either $\circ$ or $\Delta.$ In turn, if it does not admit $\circ\Rightarrow\Delta$ at that,
then only its first strong component can be $\Delta.$
Hence, it is either $TT_{\delta}$ or $\Delta\Rightarrow TT_{\delta-3}.$
This implies that the out-set of each vertex of $T$
belongs to $\Cal{T}_{3}\Rightarrow TT_{\delta-3}.$

Note that the converse $T^{-}$ of $T$ (obtained from $T$ by reversing all of its arcs)
is also a maximizer and hence, the out-set of each of its vertices is contained in
$\Cal{T}_{3}\Rightarrow TT_{\delta-3}.$
The out-set of a vertex in $T^{-}$ is the converse of the in-set
of the same vertex in $T.$ So, the in-set of each vertex of $T$ belongs to
$TT_{\delta-3}\Rightarrow\Cal{T}_{3}.$
In particular, for each $i,$
the in-set $N^{-}(i)$ contains a source.
Denote it by $j.$ Obviously, $N^{-}(j)$ is a subset of $N^{+}(i).$
Since $|N^{+}(i)|=|N^{-}(j)|=\delta,$
we have $N^{+}(i)=N^{-}(j).$
As for $\delta\ge 4,$
$\Cal{T}_{3}\Rightarrow TT_{\delta-3}\bigcap
TT_{\delta-3}\Rightarrow\Cal{T}_{3}=\{TT_{\delta}\},$
the structure of $N^{+}(i)$ and $N^{-}(j)$
described above implies that $N^{+}(i)\cong TT_{\delta}.$
Hence, $T\cong RLT_{n}.$ The proposition is proved.

\smallskip

For comparison with $(4),$ note that
$$c_{5}(RLT_{n})=\frac{(n+1)n(n-1)(n-3)(3n-11)}{480}.$$
This expression for $c_{5}(RLT_{n})$ was first obtained in [5], while
it is also presented in [28].
Let $\Delta(\circ,TT_{3},\circ)$ be the $5$-tournament obtained from $\Delta$
by replacing one of its vertices with $TT_{3}$ and
$\Delta(TT_{2},\circ,TT_{2})$ be the $5$-tournament obtained from $\Delta$ by
replacing two of its vertices with $TT_{2}$. If we add $RLT_{5}$ to them, then we obtain
a list of all strong locally transitive $5$-tournaments
(only they can be strong $5$-subtournaments of $RLT_{n}$!).
Either of the first two $5$-tournaments
admits exactly one hamiltonian cycle, while $RLT_{5}$ contains two $5$-cycles.
So,
$$c_{5}(RLT_{n})=n_{\Delta(\circ,TT_{3},\circ)}(RLT_{n})+
n_{\Delta(TT_{2},\circ,TT_{2})}(RLT_{n})+2n_{RLT_{5}}(RLT_{n}),$$
while
$$s_{5}(RLT_{n})=n_{\Delta(\circ,TT_{3},\circ)}(RLT_{n})+
n_{\Delta(TT_{2},\circ,TT_{2})}(RLT_{n})+n_{RLT_{5}}(RLT_{n}).$$
The known expressions for $s_{5}(RLT_{n})$ and $c_{5}(RLT_{n})$ imply that
$$n_{RLT_{5}}(RLT_{n})=\frac{(n+3)(n+1)n(n-1)(n-3)}{1920}.$$
Theorem 1.2 [13] allows us to suggest that the maximum of $n_{RLT_{5}}(T)$
in the class $\Cal{T}_{n}$ is attained at $RLT_{n}.$ Based on this theorem,
one can also conjecture that
$$n_{\Delta(\circ,TT_{3},\circ)}(T)\le
n_{\Delta(\circ,TT_{3},\circ)}(RLT_{n})=\frac{(n+1)n(n-1)(n-3)(n-5)}{384}.$$
However, Theorem 1.3 [13] means that for large odd $n,$ the maximum of the quantity $n_{\Delta(TT_{2},\circ,TT_{2})}(T)$
in the class $\Cal{T}_{n}$ cannot be achieved for $RLT_{n},$
while $RLT_{n}$ has the same number of copies of $\Delta(\circ,TT_{3},\circ)$
and $\Delta(TT_{2},\circ,TT_{2})$.
The results of [13] show that the problem of determining
the maximum number of copies of a given $5$-tournament in the class $\Cal{T}_{n}$
is trivial only for $TT_{5}.$ For other $5$-tournaments, it seems to be very difficult.

Denote by $\Delta\cdot\Delta$ the tournament of order $9$ obtained from
the cyclic triple $\Delta$ by replacing each of its vertices with a copy
of $\Delta.$ Lemma 1 and
formula $(3)$ allow us not only to determine all maximizers of
$s_{5}(T)$ in the class $\Cal{T}_{n}$ but also to obtain
a sharp lower bound on $s_{5}(T)$ in the class ${\Cal R}_{n}.$

\smallskip

{\bf Proposition 2.} {\sl Let $T$ be a regular tournament of (odd) order $n.$
If $n\equiv3\mod 4,$ then the
inequality}
$$\frac{n(n+1)(n-1)(n-3)(17n-59)}{3840}=s_{5}(\Cal{DR}_{n})\le s_{5}(T)$$
{\sl holds with equality iff $T$ is doubly-regular
(i.e. $T\in \Cal{DR}_{n}$) or
$T$ is an arbitrary regular $7$-tournament.
In turn,
if $n\equiv1\mod 4,$ then we have
$$\frac{n(n-1)(17n^{3}-93n^{2}+127n-243)}{3840}=s_{5}(\Cal{NDR}_{n})\le s_{5}(T)$$
with equality holding iff $T$ is nearly-doubly-regular
(i.e. $T\in \Cal{NDR}_{n}$) or $T$ is
isomorphic to $\Delta\cdot\Delta$
or one of the following two regular tournaments of order $9$
with adjacency matrices}
$$\pmatrix
0 & 1 & 1 & 1 & 1 & 0 & 0 & 0 & 0\\
0 & 0 & 1 & 1 & 1 & 1 & 0 & 0 & 0\\
0 & 0 & 0 & 1 & 0 & 1 & 0 & 1 & 1\\
0 & 0 & 0 & 0 & 1 & 0 & 1 & 1 & 1\\
0 & 0 & 1 & 0 & 0 & 0 & 1 & 1 & 1\\
1 & 0 & 0 & 1 & 1 & 0 & 1 & 0 & 0\\
1 & 1 & 1 & 0 & 0 & 0 & 0 & 1 & 0\\
1 & 1 & 0 & 0 & 0 & 1 & 0 & 0 & 1\\
1 & 1 & 0 & 0 & 0 & 1 & 1 & 0 & 0\endpmatrix \ \
\text{\sl and }\ \
\pmatrix
0 & 1 & 1 & 1 & 1 & 0 & 0 & 0 & 0\\
0 & 0 & 1 & 0 & 1 & 0 & 1 & 1 & 0\\
0 & 0 & 0 & 1 & 1 & 0 & 1 & 0 & 1\\
0 & 1 & 0 & 0 & 0 & 1 & 0 & 1 & 1\\
0 & 0 & 0 & 1 & 0 & 1 & 0 & 1 & 1\\
1 & 1 & 1 & 0 & 0 & 0 & 0 & 0 & 1\\
1 & 0 & 0 & 1 & 1 & 1 & 0 & 0 & 0\\
1 & 0 & 1 & 0 & 0 & 1 & 1 & 0 & 0\\
1 & 1 & 0 & 0 & 0 & 0 & 1 & 1 & 0\endpmatrix.$$

\smallskip
{\bf Proof.} For a regular tournament $T$ with semi-degree $\delta$
(and hence, of order $n=2\delta+1$) and $m=5,$ formula $(3)$ takes the following form
$$s_{5}(T)=
\binom{n}{5}-2n\binom{\delta}{4}
+\sum\limits_{i=1}^{n}n_{\Cal{T}_{3}\Rightarrow\circ}\Bigl(N^{+}(i)\Bigr).\eqno(6)$$
This formula, equality $(2),$ and Lemma 1 imply that for odd $\delta\ge 1,$ we have
$$s_{5}(T)\ge \binom{n}{5}-2n\binom{\delta}{4}+
n\delta\binom{\frac{\delta-1}{2}}{3}
=\frac{\delta(\delta-1)(\delta+1)(2\delta+1)(17\delta-21)}{240}.$$
For odd $\delta\ge 5$,
the equality holds if and only if for each $i=1,...,n,$ the out-set
$N^{+}(i)$ is a regular tournament of order $\delta.$
In turn, if $\delta=3,$ then any $N^{+}(i)$ satisfies
$n_{\Cal{T}_{3}\Rightarrow\circ}\Bigl(N^{+}(i)\Bigr)=0.$
So, for any regular tournament $T$ of order $7,$ the equality holds.
The case $\delta=1$ is trivial as $\Delta$ is the unique regular tournament of order $3$
and it is also the unique element of $\Cal{DR}_{3}.$

Formula $(6),$ equality $(2),$ and Lemma 1 imply that
for even $\delta\ge 2$, we have $$s_{5}(T)\ge\binom{n}{5}-2n\binom{\delta}{4}+
n\frac{\delta}{2}\Bigl(\binom{\frac{\delta}{2}}{3}+\binom{\frac{\delta}{2}-1}{3}\Bigr)=
\frac{\delta(2\delta+1)(17\delta^{3}-21\delta^{2}-2\delta-24)}{240}.$$
For even $\delta\ge 6$,
the equality holds if and only if for each $i=1,...,n,$ the out-set
$N^{+}(i)$ is a near regular tournament of order $\delta.$
For $\delta=4,$ this
lower bound is attained if and only if for each $i,$
the out-set $N^{+}(i)$ is not contained in $\Cal{T}_{3}\Rightarrow\circ,$ i.e.
it is either $ST_{4}$ or $\circ\Rightarrow\Delta.$
It is not difficult to check that the out-set
of each vertex in the composition $\Delta\cdot\Delta$
induces $\circ\Rightarrow \Delta.$ It is also shown in [33]
that there are exactly two nearly-doubly-regular tournaments of order $9$
(the out-set of each of their vertices induces $ST_{4}$).
Can one present other examples of $T\in {\Cal R}_{9}$ the out-sets
of whose vertices induce $ST_{4}$ or $\circ\Rightarrow \Delta$?

Note that the score-lists of $ST_{4}$ and $\circ\Rightarrow\Delta$ are
$(1,1,2,2)$ and $(1,1,1,3),$ respectively.
The score-lists of the remaining tournaments of order $4$ contain $0.$
Non-diagonal entries of the matrix $AA^{\top},$ where $A$ is the adjacency
matrix and $A^{\top}$ is its transpose, form the score-lists of
the out-sets of the vertices.
So, we have to find regular tournaments $T$ of order $9$
for which $AA^{\top}>0,$ i.e. the matrix $AA^{\top}$ has only positive entries.
According to [3], there exist exactly
15 (non-isomorphic) regular tournaments of order $9.$
Our exhaustive manual
consideration of the list of all elements of ${\Cal R}_{9}$ presented therein shows that
there are another two regular tournaments of order $9$ with $AA^{\top}>0$.
Their adjacency matrices were given in the statement of the proposition.
Note that for either of them,
there exists a vertex whose out-set induces $ST_{4}$ and
there exists a vertex whose out-set induces $\circ\Rightarrow \Delta.$
So, in some sense, they are placed between two elements of $\Cal{NDR}_{9}$ and
$\Delta \cdot \Delta.$ Finally, the case $\delta=2$ is trivial
because the unique regular tournament
of order $5$ is nearly-doubly-regular (it is isomorphic to $RLT_{5}$) and hence,
the statement also holds for $n=5.$
The proposition is proved.

As we have seen in the introduction, formula $(3)$ and Lemma 1 imply that
for $T\in \Cal{R}_{n},$ the inequality $s_{4}\bigl(\Cal{(N)DR}_{n}\bigr) \le s_{4}(T)$
or, the same, $n_{ST_{4}}\bigl(\Cal{(N)DR}_{n}\bigr) \le n_{ST_{4}}(T)$
holds with equality if and only if $T\in \Cal{(N)DR}_{n}.$ The same also takes place
for the quantity $n_{TT_{4}}(T).$ For $\circ\Rightarrow\Delta$ or
$\circ\Leftarrow\Delta,$ the minimum number of copies of either of them
is attained at $RLT_{n}$ and equals $0.$ The same also holds for all tournaments
of order $5$ with the exception of $TT_{5},$ $RLT_{5},$ $\Delta(\circ,TT_{3},\circ),$
and $\Delta(TT_{2},\circ,TT_{2})$ because only they are locally transitive and hence,
only they can be 5-subtournaments of $RLT_{n}.$ For each of them,
the problem of determining the minimum number of copies in the class $\Cal{R}_{n}$
is difficult. It is so even for $TT_{5}.$ Indeed, the identity
$$n_{TT_{5}}(T)=\sum\limits_{i=1}^{n}n_{TT_{4}}\Bigl(N^{+}(i)\Bigr)=
\sum\limits_{i=1}^{n}n_{TT_{4}}\Bigl(N^{-}(i)\Bigr)$$
and the lower bound $n_{TT_{4}}\bigl(\Cal{(N)DR}_{n}\bigr) \le n_{TT_{4}}(T)$
allow us to get a lower bound on $n_{TT_{5}}(T).$
(For $n\equiv 7\mod 8,$ the corresponding expression is given in the last
section of [20].)
For $n\ge 11,$ it is attained only at a {\sl regular locally doubly-regular tournament}
(when $n\equiv 7\mod 8$), i.e. the out-set (in-set) of each vertex is doubly regular,
or a {\sl regular locally nearly-doubly-regular tournament} (when $n\equiv 3\mod 8$), i.e. the out-set (in-set)
of each vertex is nearly doubly regular. Denote by $\Cal{RLDR}_{n}$ and $\Cal{RLNDR}_{n}$
the corresponding subfamilies of $\Cal{R}_{n}.$ Are they non-empty?

For $n=p^{k},$ where $p$ is
a prime that is congruent to $3$ modulo $4$ and $k$ is an odd positive integer,
define the {\sl quadratic residue tournament} $QR_{n}$ as the tournament such that
its vertex-set is the Galois field $GF(n)$ and $i\to j$ if and only if $j-i$
is a non-zero square in $GF(n).$ For each such $n,$ the tournament $QR_{n}$
has high regularity properties. In particular, it
is doubly regular.
One can check that $\Cal{RLDR}_{7}=\{QR_{7}\},$ while $\Cal{RLDR}_{15}=\emptyset$
and for $n=23,31,$ and $47,$ the tournament $QR_{n}$ does not belong to
the set $\Cal{RLDR}_{n}.$
Moreover, $\Cal{RLNDR}_{11}=\{QR_{11}\},$ $\Cal{RLNDR}_{19}=\{QR_{19}\},$
and $\Cal{RLNDR}_{27}$ contains $QR_{27},$ while $QR_{43}$ is not contained
in $\Cal{RLNDR}_{43}.$ We conjecture that a regular locally doubly-regular tournament
of order $n$ does not exist for any $n>7$ (and so, the Moon lower bound on $n_{TT_{5}}(T)$
is not sharp for $n>7$)
and for sufficiently large $n$, the set
$\Cal{RLNDR}_{n}$ is also empty. For these values of $n,$ we have no good candidates
for minimizers of $n_{TT_{5}}(T)$. At the moment, we can only state that
any sequence of minimizers of growing order $n$ must be
quasi-random (see, for instance, [11]) or, the same, asymptotically
doubly regular (the corresponding definition will be given in Section 4 below).

Proposition 2 and the main result of [4] imply that for $m=4$ or $5$ and
$T\in {\Cal R}_{n},$ we have
$s_{m}\bigl(\Cal{(N)DR}_{n}\bigr) \le s_{m}(T)\le s_{m}(RLT_{n}).$
On the other hand,
it is shown in [28] (see also Proposition A.1 in Appendix A) that
$$c_{5}(T)+2c_{4}(T)=\frac{n(n-1)(n+1)(n-3)(n+3)}{160}\eqno(7)$$
and hence,
$c_{5}(RLT_{n})\le c_{5}(T)\le c_{5}\bigl(\Cal{(N)DR}_{n}\bigr).$
So, while $c_{m}(T)$ and $s_{m}(T)$ coincide for $m=3$ and $m=4,$
they demonstrate quite different
behaviour for the next value $m=5.$
In the next section, we prove that the inequality
$c_{5}(T)\le c_{5}(\Cal{DR}_{n})$ holds for each $T\in {\Cal T}_{n},$
where $n$ is odd.

\bigskip
\centerline{\bf \S 3. An upper bound on $\bold{c_{5}(T)}$ in the class $\bold{{\Cal T}_{n}},$
where $\bold{n}$ is odd}

\medskip

For $m=4,$ formula $(3)$ can be rewritten as
$$c_{4}(T)=\binom{n}{4}-\sum\limits_{i=1}^{n}\binom{\delta_{i}^{-}}{3}
-\sum\limits_{i=1}^{n}c_{3}\Bigl(N^{+}(i)\Bigr)$$
or
$$c_{4}(T)=\binom{n}{4}-\sum\limits_{i=1}^{n}\binom{\delta_{i}^{+}}{3}
-\sum\limits_{i=1}^{n}c_{3}\Bigl(N^{-}(i)\Bigr).$$
Replacing all $\delta_{i}^{+}$ and $\delta_{i}^{-}$ by $\frac{n-1}{2},$
summing both identities and then dividing the sum obtained by $2$ yield
a formula for the number of $4$-cycles in $T\in \Cal{R}_{n}$:
$$c_{4}(T)=\frac{(n+1)n(n-1)(n-3)}{48}-
\frac{1}{2}\sum\limits_{i=1}^{n}c_{3}
\Bigl(N^{+}(i)\Bigr)
-\frac{1}{2}\sum\limits_{i=1}^{n}c_{3}
\Bigl(N^{-}(i)\Bigr).$$

In turn, this formula and identity $(7)$
imply that for a regular tournament $T$ of order $n,$ we have
$$c_{5}(T)=\frac{n(n-1)(n+1)(n-3)(3n-11)}{480}+
\sum\limits_{i=1}^{n}c_{3}\bigl(N^{+}(i)\bigr)+
\sum\limits_{i=1}^{n}c_{3}\bigl(N^{-}(i)\bigr).$$
We see that two regular $n$-tournaments with the same collection of
the out-sets and in-sets of vertices have the same number of $5$-cycles.
However, it is not so even for near regular locally
transitive tournaments of given even order $n.$  In particular,
$\Delta\cdot TT_{2}$ has $6$ cycles of length $5,$ while $RLT_{7}-\circ$
(it is obtained by replacing one vertex of $RLT_{5}$ with $TT_{2}$) admits
$8$ cycles of length $5.$ This is a direct consequence of the following three
facts: $(a)$ any near regular tournament of order $6$ is $2$-strong;
$(b)$ $c_{5}(RLT_{5})=2;$ $(c)$ any other strong locally
transitive tournament of order $5$ admits exactly one hamiltonian cycle.

As one can see, generally,
the quantity $c_{5}(T)$ is not uniquely
determined by the out-sets of the vertices.  In [16], a formula for $c_{5}(T)$
was obtained  in terms of the {\sl intersection numbers} (or, the {\sl edge-scores}).
\footnote[1]{As the authors of [16] remark, the combinatorial sense of their formula
is not well understood. We would also like to note that while
the known expressions for $c_{5}(T)$ in the classes $\Cal{R}_{n}$ and $\Cal{NR}_{n}$
contain only the numbers of $3$-cycles
in some subtournaments, at the moment, we are unable to present such a formula
for $c_{5}(T)$ in the general case (recall that it exists for $c_{4}(T)$).}
For any two vertices $i$ and $j$, they are defined as follows
$$\delta_{ij}^{++}=|N^{+}(i)\cap N^{+}(j)|, \ \ \ \delta_{ij}^{+-}=|N^{+}(i)\cap N^{-}(j)|,$$
and
$$\delta_{ij}^{--}=|N^{-}(i)\cap N^{-}(j)|,\ \ \ \delta_{ij}^{-+}=|N^{-}(i)\cap N^{+}(j)|.$$

Formulas $(2)$ and $(3)$ imply that
$$s_{5}(T)=\binom{n}{5}-\sum\limits_{i=1}^{n}
\binom{\delta^{-}_{i}}{4}-
\sum\limits_{i=1}^{n}\binom{\delta^{+}_{i}}{4}
+\sum\limits_{(i,j)\in {\Cal A}(T)}\binom{\delta^{+-}_{ij}}{3},$$
where ${\Cal A}(T)$ is the arc-set of $T$.
The formula for $c_{5}(T)$ can be written in our notation as follows:
$$8c_{5}(T)=6\binom n5+\sum\limits_{(i,j)\in {\Cal A}(T)}
f(\delta_{ij}^{++},\delta_{ij}^{--},\delta_{ij}^{-+},\delta_{ij}^{+-}),$$
where
$$f(\delta_{ij}^{++},\delta_{ij}^{--},\delta_{ij}^{-+},\delta_{ij}^{+-})=
-(\delta_{ij}^{+-}+\delta_{ij}^{-+})
(\delta_{ij}^{++}-\delta_{ij}^{--})^{2}-
(\delta_{ij}^{++}+\delta_{ij}^{--})
(\delta_{ij}^{+-}-\delta_{ij}^{-+})^{2}+$$
$$2(\delta_{ij}^{++}+\delta_{ij}^{--})(\delta_{ij}^{+-}+\delta_{ij}^{-+}).\eqno(8)$$
Note that
$$\delta_{ij}^{++}+\delta_{ij}^{--}+\delta_{ij}^{+-}+\delta_{ij}^{-+}=n-2.\eqno(9)$$
Based on the formula for $c_{5}(T)$ and $(9)$, the authors of [16] showed that
$$c_{5}(T)\le \frac{3}{4}\binom{n}{5}+
\frac{1}{4}\binom{n}{2}\Bigl(\frac{n-2}{2}\Bigr)^{2}=
\frac{n(n-1)(n-2)(n^{2}-2n+2)}{160}.$$
A more detailed analysis of $(8)$ allows us to obtain the following
proposition.

\smallskip

{\bf Theorem 1.} {\sl For a tournament $T$ of odd order $n\ge 5,$ the
inequality}
$$c_{5}(T)\le \frac{(n+1)n(n-1)(n-2)(n-3)}{160}\eqno(10)$$
{\sl holds with equality iff $T$ is doubly-regular (i.e. $T\in \Cal{DR}_{n}$).}

\smallskip
{\bf Proof.}
Let us first give an upper bound on the value of the function
$f(\delta_{ij}^{++},\delta_{ij}^{--},$ $\delta_{ij}^{-+},\delta_{ij}^{+-})$ in $(8)$
for a given arc $(i,j)$ of $T.$
Since $n-2$ is odd, equality $(9)$ implies that
at least one of inequalities  $\delta_{ij}^{-+}\neq \delta_{ij}^{+-}$
and $\delta_{ij}^{++}\neq \delta_{ij}^{--}$ holds.

Assume first that $\delta_{ij}^{-+}\neq \delta_{ij}^{+-}$
(and hence, $\delta_{ij}^{+-}+\delta_{ij}^{-+}\ge 1$).
Then $(8)$ implies that
$$f(\delta_{ij}^{++},\delta_{ij}^{--},\delta_{ij}^{-+},\delta_{ij}^{+-})
\le
-\bigl(\delta_{ij}^{++}+\delta_{ij}^{--}\bigr)+2\bigl(\delta_{ij}^{++}+\delta_{ij}^{--}\bigr)
\bigl(\delta_{ij}^{+-}+\delta_{ij}^{-+}\bigr)$$
with equality holding iff
$\delta_{ij}^{++}=\delta_{ij}^{--}\ge 1$ and $\delta_{ij}^{+-}=\delta_{ij}^{-+}\pm 1$
or $\delta_{ij}^{++}=\delta_{ij}^{--}=0$.
Let $\Delta_{ij}=\delta_{ij}^{++}+\delta_{ij}^{--}.$ By $(9),$ we have
$\delta_{ij}^{+-}+\delta_{ij}^{-+}=n-2-\Delta_{ij}.$
So, the above inequality
can be rewritten as
$$f(\delta_{ij}^{++},\delta_{ij}^{--},\delta_{ij}^{-+},\delta_{ij}^{+-})\le
-\Delta_{ij}+2\Delta_{ij}\bigl(n-2-\Delta_{ij}\bigr).$$
Note that
$$-\Delta_{ij}+2\Delta_{ij}\bigl(n-2-\Delta_{ij}\bigr)=
\frac{(n-3)(n-2)}{2}-2\Bigl(\Delta_{ij}-\frac{n-3}{2}\Bigr)
\Bigl(\Delta_{ij}-\frac{n-2}{2}\Bigr).$$
For integers $\Delta_{ij}$ and $n,$ the subtrahend in the right-hand side of the equality
is always non-negative. It is equal to zero
iff $\Delta_{ij}=\frac{n-3}{2}$ or $\Delta_{ij}=\frac{n-2}{2}.$
The first number is strictly greater than $0$ if $n\ge 5$
and the last number is not an integer if $n$ is odd.
Hence, for the case $\delta_{ij}^{-+}\neq \delta_{ij}^{+-},$ we have
$$f(\delta_{ij}^{++},\delta_{ij}^{--},\delta_{ij}^{-+},\delta_{ij}^{+-})\le
\frac{(n-3)(n-2)}{2}\eqno(11)$$
with equality holding iff either
$$\delta_{ij}^{++}=\delta_{ij}^{--}=\delta_{ij}^{+-}=\frac{n-3}{4}\
\text{ and }\ \delta_{ij}^{-+}=\frac{n+1}{4}\ \ \ \bigl(\text{case}\
\delta_{ij}^{-+}=\delta_{ij}^{+-}+1\bigr)\eqno(a)$$
or
$$\delta_{ij}^{++}=\delta_{ij}^{--}=\delta_{ij}^{-+}=\frac{n-3}{4}\
\text{ and }\ \delta_{ij}^{+-}=\frac{n+1}{4}\ \ \ \bigl(\text{case}\
\delta_{ij}^{-+}=\delta_{ij}^{+-}-1\bigr).\eqno(b)$$

Note that the value of $f(\delta_{ij}^{++},\delta_{ij}^{--},\delta_{ij}^{-+},\delta_{ij}^{+-})$
remains the same if one interchanges the pairs $(\delta_{ij}^{++},\delta_{ij}^{--})$
and $(\delta_{ij}^{+-},\delta_{ij}^{-+}).$
Thus, the arguments given above imply that
for the case $\delta_{ij}^{++}\neq \delta_{ij}^{--},$ inequality $(11)$
also holds with equality
iff  either
$$\delta_{ij}^{+-}=\delta_{ij}^{-+}=\delta_{ij}^{--}=\frac{n-3}{4}
\ \text{ and }\ \delta_{ij}^{++}=\frac{n+1}{4}\ \ \
\bigl(\text{case}\ \delta_{ij}^{++}=\delta_{ij}^{--}+1\bigr)
\eqno(c)$$
or
$$\delta_{ij}^{+-}=\delta_{ij}^{-+}=\delta_{ij}^{++}=\frac{n-3}{4}
\ \text{ and }\ \delta_{ij}^{--}=\frac{n+1}{4}\ \ \
\bigl(\text{case}\ \delta_{ij}^{++}=\delta_{ij}^{--}-1\bigr)
.\eqno(d)$$

Inequality $(11)$ taken together with $(8)$ implies that
$$c_{5}(T)\le
\frac{3}{4}\binom n5+\frac{1}{8}\sum\limits_{(i,j)\in {\Cal A}(T)}
\frac{(n-3)(n-2)}{2}=$$
$$=\frac{n(n-1)(n-2)(n-3)(n-4)}{160}+\frac{n(n-1)(n-2)(n-3)}{32}=$$
$$=\frac{n(n-1)(n-2)(n-3)(n+1)}{160}.$$

Assume that this upper bound is attained.
Then for each arc $(i,j)\in {\Cal A}(T),$ the equality in $(11)$
holds and hence, exactly one of cases $(a)-(d)$ is realized.
In all cases, we have
$\delta_{ij}^{++}\ge \frac{n-3}{4}.$  Note that $\delta_{ij}^{++}$ is the out-degree
of $j$ in the subtournament induced by $N^{+}(i).$ For a tournament
with minimum out-degree $\delta^{+}_{min}$ (maximum out-degree
$\delta^{+}_{max}$), its order is at least $2\delta^{+}_{min}+1$
(at most $2\delta^{+}_{max}+1$) with equality holding iff it is
regular.
Hence, for each vertex
$i$ of $T,$ the inequality $\delta_{i}^{+}\ge \frac{n-1}{2}$ holds. As the order
of $T$ is $n,$ in fact, we have $\delta_{i}^{+}=\frac{n-1}{2}$ for each $i$ and hence,
$\delta_{ij}^{++}=\frac{n-3}{4}$ for each $(i,j)\in {\Cal A}(T),$
i.e. $T$ is doubly regular.

Note that if $i\to j,$ then the intersection numbers
are related as follows:
$$\delta_{i}^{+}=1+\delta_{ij}^{++}+\delta_{ij}^{+-},\ \ \
\delta_{i}^{-}=\delta_{ij}^{--}+\delta_{ij}^{-+}$$
and
$$\delta_{j}^{+}=\delta_{ij}^{++}+\delta_{ij}^{-+},\ \ \
\delta_{j}^{-}=1+\delta_{ij}^{--}+\delta_{ij}^{+-}.$$
From this it follows that if
$$\delta_{i}^{+}=\delta_{i}^{-}=\delta_{j}^{+}=\delta_{j}^{-}=\frac{n-1}{2}\ \
\text{ and } \ \delta_{ij}^{++}=\frac{n-3}{4},$$
then
$$\delta_{ij}^{++}=\delta_{ij}^{--}=\delta_{ij}^{+-}=\frac{n-3}{4}\
\text{ and }\ \delta_{ij}^{-+}=\frac{n+1}{4}.$$
This means that
for each arc $(i,j)$ of $T\in \Cal{DR}_{n},$ case $(a)$ always holds.
Thus, the upper bound is achieved if and only if $T\in \Cal{DR}_{n}.$
The theorem is proved.

\smallskip

Unfortunately, for $n\equiv 1\mod 4,$ upper bound $(10)$ is not sharp.
It was shown in [28] that for this case, in the class
${\Cal R}_{n}$, an upper bound
$$c_{5}(T)\le \frac{n(n-1)(n^{3}-4n^{2}+n-14)}{160}
=c_{5}(\Cal{NDR}_{n})\eqno(12)$$
holds with equality if and only if $T\in \Cal{NDR}_{n}.$
However, inequality $(12)$ does not hold in the class $\Cal{T}_{5}$
because the number of hamiltonian cycles in
the tournament $\Delta(\circ,\Delta,\circ)$ obtained from
the cyclic triple $\Delta$ by replacing one of its vertices with a copy of it
is equal to the number of hamiltonian paths in $\Delta$ and so,
$c_{5}\bigl(\Delta(\circ,\Delta,\circ)\bigr)=3,$ while
according to $(12),$ we have $c_{5}(\Cal{NDR}_{5})=2.$

Nevertheless, our computer search shows that the maximum numbers of $5$-cycles
in ${\Cal T}_{n}$ and ${\Cal R}_{n}$ coincide for $n=9.$
We believe that the Komarov-Mackey formula for $c_{5}(T)$ and
arguments which are similar to those used in the proof of Theorem 1
will allow us to prove that
for $n\equiv1\mod4$ and $n>5,$
upper bound $(12)$ also holds in the class ${\Cal T}_{n},$
but we will not try to do this here because without any doubt,
a possible proof for this case contains much more routine calculations
and hence, should be considered in a separate paper.

We see that for $m=5,$ the problem of determining the maximum number of $m$-cycles
in the class ${\Cal T}_{n}$ is not so simple as that in the case of $m=3$ or $m=4.$
Note that the length $m=5$ is critical
for many combinatorial problems involving $m$-cycles. In particular, the question of decomposition
of complete tripartite graphs into cycles of odd length $m$
also becomes non-trivial for $m=5$ (see [19]).

\newpage
\centerline{\bf \S 4. Concluding remarks on $\bold{s_{m}(T)}$
and $\bold{c_{m}(T)}$ in the case of arbitrary $\bold{m}$}
\medskip

For arbitrary $m\ge 3,$ the right-hand side of $(3)$ is equal to
the number $w_{m}(T)$ of subtournaments of order $m$ containing neither sink nor source in $T.$
By Lemma 1, for each $m\ge 3$ and odd $n$ with $\exists$-property,
the minimum of $w_{m}(T)$ in the class $\Cal{R}_{n}$ is attained at $\Cal{(N)DR}_{n}.$
For $m=3,4,$ and $5,$ we have $w_{m}(T)=s_{m}(T).$ However,
it is not so for $m=6.$ For getting an expression for $s_{m}(T)$ in this case, we have
to subtract $n_{\Delta\Rightarrow\Delta}(T)$ from the right-hand side of $(3).$
To obtain an explicit formula for $n_{\Delta\Rightarrow\Delta}(T),$ we need to consider
the intersection numbers of higher order, namely,
$\delta_{ikj}^{\pm\pm\pm}=|N^{\pm}(i)\cap N^{\pm}(k)\cap N^{\pm}(j)|.$
As a consequence, a possible expression for $s_{m}(T)$ becomes more complicated
for $m\ge 6.$ However, we think that the necessary corrections to the right-hand side of $(3)$
do not essentially change the situation
and conjecture that the minimum of $s_{m}(T)$ in the class $\Cal{R}_{n}$
is also achieved for some element of $\Cal{(N)DR}_{n}.$

Note that the connectivity number of a regular tournament of order $n$ is at least
$\lceil \frac{n}{3}\rceil$ (see Lemma 4.1 [34]).
Hence, for any $n$ that is sufficiently close to $m,$
any two regular $n$-tournaments have the same number of strong $m$-subtournaments,
namely, $\binom{n}{m},$ and so, in this case,
we can write $s_{m}(\Cal{R}_{n})=\binom{n}{m}$ (when $n\approx m$).
However, it is not so if $n$ is large enough.
We suggest that for arbitrary $m\ge 4$ (not only for $m=4$ and $m=5$) and sufficiently
large odd $n$ with $\exists$-property,
any minimizer of $s_{m}(T)$ in the class $\Cal{R}_{n}$
is contained in $\Cal{(N)DR}_{n},$ while the maximum of $s_{m}(T)$ in the class
$\Cal{T}_{n}$ is attained only at $RLT_{n}.$ The latter can be also conjectured
for $n_{RLT_{m}}(T)$ if $m$ is odd (see Conjecture 5.5 in [12] and Conjecture 3
in [6]). The known expressions for this quantity in the case of $m=3$
and $m=5$ allow us to suggest that
$$n_{RLT_{m}}(RLT_{n})=$$
$$\frac{(n+m-2)(n+m-4)...(n+1)n(n-1)...(n-m+4)(n-m+2)}{2^{m-1}m!}.$$
In our further papers, we will try to confirm this conjecture.

Note that the trace $tr_{m}(T)$ of the $m$th power of
the adjacency matrix of $T$ is equal to the number of closed $m$-walks on $T$
and hence, $mc_{m}(T)\le tr_{m}(T).$
For $m=3,4,$ and $5,$ any such walk can be obtained
by a shift along some $m$-cycle. Hence, for these values of $m,$ we have
$mc_{m}(T)=tr_{m}(T)$ for each $T.$  However, it is not so for $m=6$ because
repeating a closed $3$-walk ($3$-cycle) provides a closed $6$-walk.
Nevertheless, for arbitrary $m$, there exists $C_{m}>0$ not depending on $n$
such that
$$tr_{m}(T)-mc_{m}(T)< C_{m}n^{m-1}.\eqno(13)$$
Inequality $(13)$ allows us
to study the asymptotical properties of $c_{m}(T)$ in different subfamilies of
$\Cal{T}_{n}$ as $n\to \infty$ with the use of
the known spectral properties of tournament matrices.

Let $R_{n}^{(m)}$ be a regular tournament of order $n$
which maximizes the number of cycles of length $m\ge 4$ in the class $\Cal{R}_{n}.$
It is shown in [32] that $R_{n}^{(m)}$ is {\sl asymptotically doubly regular}
(i.e. $n_{\circ\Rightarrow \Delta}(R_{n}^{(m)})=
n_{\circ\Rightarrow \Delta}(\Cal{DR}_{n})+o(n^{4});$ here,
$\circ\Rightarrow \Delta$ can be replaced by any element of $\Cal{T}_{4}$)
or, the same, {\sl quasi-random} if and only if $m$ is not a multiple of $4.$
Our Theorem 1 means that $c_{5}(T)\le c_{5}(\Cal{DR}_{n})$ for each $T\in \Cal{T}_{n},$
where $n$ is odd. In turn, Theorem 4 [28] shows that
$c_{6}(T)\le c_{6}(\Cal{DR}_{n})$ for each $T\in \Cal{R}_{n},$
where $n\ge 7.$ Finally, in [29], many serious arguments
are given for supporting the conjecture that
$c_{7}(T)\le c_{7}(\Cal{DR}_{n})$ for each $T\in \Cal{R}_{n},$
where $n\ge 7.$
Note that at the moment, $(5,5)$ is the only known pair $(m,n)$ for which
the maxima of $c_{m}(T)$ in ${\Cal T}_{n}$ and ${\Cal R}_{n}$
are distinct.
All these facts and also the results of [14] allow us to suggest that
for each $m\equiv 1,2,3\mod 4$ and sufficiently large
odd $n$ with $\exists$-property, we have
$$\max\{c_{m}(T):T\in \Cal{T}_{n}\}=\max\{c_{m}(T):T\in \Cal{(N)DR}_{n}\}.\eqno(14)$$

For these values of $m,$
the maximum of $c_{m}(T)$ in the class ${\Cal R}_{n}$
(and even in the class ${\Cal T}_{n}$ as one of the main theorems of [14] states) is
asymptotically the same as the expected number $E_{n}c_{m}$ of $m$-cycles in a random $n$-tournament.
According to G. Korvin (see [17]), for $m\ge 3,$ we have
$E_{n}c_{m}=\bigl(n\bigr)_{m}/\bigl(m2^{m}\bigr),$
where $\bigl(n\bigr)_{m}=n(n-1)...(n-m+1).$
It is shown in [24] and [27]
that for $m=5,$ the quantity $c_{m}(\Cal{DR}_{n})$ is equal to
$E_{n+1}c_{m}$ (see also Theorem 1).
The same also holds for $m=3$ (see [15]).
The known expressions for
$c_{m}(\Cal{DR}_{n})$ in the case of $m=6,7,8,$ and $9$
obtained in [28], [29], [30], and [31], respectively,
imply that for $6\le m\le 9$ and $n\ge m,$ we have
$$\bigl(n\bigr)_{m}/\bigl(m2^{m}\bigr)=E_{n}c_{m}< c_{m}(\Cal{DR}_{n})< E_{n+1}c_{m}=\bigl(n+1\bigr)_{m}/\bigl(m2^{m}\bigr).$$
Based on this fact, it is natural to suggest that for $m>5$ and odd $n\ge m,$
the strict inequality
$$\max\{c_{m}(T):T\in \Cal{(N)DR}_{n}\}< E_{n+1}c_{m}=\bigl(n+1\bigr)_{m}/\bigl(m2^{m}\bigr)\eqno(15)$$
always holds.
So, if both of our conjectures $(14)$ and $(15)$ are true, then
in the case of $m\equiv 1,2,3 \mod4$ and odd $n\ge m,$ for
an $n$-tournament $T,$ the
inequality $c_{m}(T)\le \bigl(n+1\bigr)_{m}/\bigl(m2^{m}\bigr)$ holds with equality
iff $m=3$ and $T$ is regular or $m=5$ and $T$ is doubly regular. This cannot
be true for $m\equiv0\mod 4$ because according to [28] and [29],
$c_{m}\bigl(RLT_{n}\bigr)=
\frac{1+(-1)^{\frac{m}{2}}\beta(m)}{m2^{m}}n^{m}+O(n^{m-1}),$
where $\beta(m)$ is the (positive) coefficient of $z^{m-1}$ in the
Maclaurin expansion of the trigonometric function $\tan z.$
According to our conjecture, for such $m$ and sufficiently large odd $n,$
the maximum of $c_{m}(T)$ in the class $\Cal{T}_{n}$ is attained at $RLT_{n}$
as in the case of $s_{m}(T).$ The results of our paper [30] show that the
condition "$n$ should be large enough" is essential even for $m=8.$ For this
case, $n$ should be strictly greater than $37$. (Recall that for $m=4,$
the inequality $c_{m}(T)\le c_{m}(RLT_{n})$ holds for each $T\in \Cal{T}_{n},$
where $n$ is an arbitrary odd number.)

\bigskip
{\bf Acknowledgments}
\smallskip

The author thanks the referees for their useful suggestions which led
to improvements in the text of the paper.
He is also grateful to N. Komarov and J. Mackey for sending him
their joint paper [16].

\bigskip
\centerline{\bf Appendix A. Identity relating $\bold{c_{4}(T)}$ and $\bold{c_{5}(T)}$}
\smallskip
\centerline{\bf in the case of a regular n-tournament $\bold{T}$}
\smallskip

We have already used an identity including
$c_{4}(T),$ $c_{5}(T),$ and $n$ for $T\in {\Cal R}_{n}$
above several times.
It was (first) obtained in [28].
In this paper, it is proved with the use of purely matrix methods.
More precisely, it is deduced from a matrix identity whose
proof essentially uses the fact that for a regular tournament,
its adjacency matrix $A$ commutes with the all ones matrix $J.$
One can also prove it with the use of the Komarov-Mackey formula for $c_{5}(T).$
In the present paper, we give a proof based on the spectral properties of $T.$
This proof is straightforward and uses no preliminary results.
By this reason, we present it here.

\smallskip

{\bf Proposition A.1 [28].} {\sl For a regular tournament $T$ of (odd) order $n,$ we have}
$$c_{5}(T)+2c_{4}(T)=\frac{n(n-1)(n+1)(n-3)(n+3)}{160}.$$

{\bf  Proof.} Let $tr_{m}(T)$ be the trace of the $m$th power of
the adjacency matrix $A$ of $T.$ Since for $m=3,4,$ and $5,$ we have
$mc_{m}(T)=tr_{m}(T),$ it suffices to relate $tr_{5}(T),$ $tr_{4}(T),$ and $n.$
The Perron root (spectral radius)
of $A$ equals $\frac{n-1}{2}$ and
its algebraic multiplicity is equal to $1.$
According to [7], all the other
(non-Perron) eigenvalues $\lambda_{1},...,\lambda_{n-1}$
lie on the line $\Re(\lambda)=-\frac{1}{2}$.
Denote by $\rho_{j}$ the imaginary part of $\lambda_{j},$ where $j=1,...,n-1.$
The binomial formula for exponent $4$ implies that
$$\lambda_{j}^{4}=\Bigl(-\frac{1}{2}+i\rho_{j}\Bigr)^{4}=
\frac{1}{16}-\frac{1}{2}i\rho_{j}-\frac{3}{2}\rho_{j}^{2}+2i\rho_{j}^{3}
+\rho_{j}^{4}.$$
In turn, the binomial formula for exponent $5$ means that
$$\lambda_{j}^{5}=\Bigl(-\frac{1}{2}+i\rho_{j}\Bigr)^{5}=
-\frac{1}{32}+\frac{5}{16}i\rho_{j}+\frac{5}{4}\rho_{j}^{2}
-\frac{5}{2}i\rho_{j}^{3}-\frac{5}{2}\rho_{j}^{4}+i\rho_{j}^{5}.$$
So,
$$\lambda_{j}^{5}+\frac{5}{2}\lambda_{j}^{4}=\frac{1}{8}-\frac{5}{2}\rho_{j}^{2}+
i\Bigl(-\frac{15}{16}\rho_{j}+\frac{5}{2}\rho_{j}^{3}+\rho_{j}^{5}\Bigr).\eqno(A.1)$$

As $A$ is a matrix with real entries, for each odd $m,$
the $m$th moment $\sum\limits_{j=1}^{n-1}\rho_{j}^{m}$ is equal to zero.
Moreover, according to [8], we have
\footnote[2]{The equality presented below is a simple consequence
of the binomial formula for exponent $2$ and the evident equality
$tr_{2}(T)=0$.}
$$\sum\limits_{j=1}^{n-1}\rho_{j}^{2}=\frac{n(n-1)}{4}.$$
Thus, summing over $j$ from $1$ to $n-1$ in $(A.1)$
and adding the terms associated with the Perron root yield
$$tr_{5}(T)+\frac{5}{2}tr_{4}(T)=
\frac{(n-1)^{5}}{32}+5\frac{(n-1)^{4}}{32}+\frac{n-1}{8}-\frac{5}{8}n(n-1)=$$
$$\frac{(n-1)}{32}\Bigl\{(n-1)^{4}+5(n-1)^{3}+4-20n\Bigr\}=
\frac{(n-1)}{32}\Bigl\{n^{4}+n^{3}-9n^{2}-9n\Bigr\}=$$
$$\frac{(n-1)(n^{2}+n)(n^{2}-9)}{32}=\frac{n(n-1)(n+1)(n-3)(n+3)}{32}.$$
Recalling that $tr_{5}(T)=5c_{5}(T)$ and $tr_{4}(T)=4c_{4}(T)$
completes the proof.

\newpage
\centerline{\bf References}
\smallskip

[1] B. Alspach and C. Tabib, A note on the number of $4$-circuits
in a tournament, Ann. Discrete Math. 12 (1982), 13-19.

\smallskip

[2] A. Asti\'e-Vidal, V. Dugat, Autonomous parts and decomposition of
regular tournaments, Discrete Math. 111 (1993), 27-36.

\smallskip

[3] A. Asti\'e-Vidal, V. Dugat, Z. Tuza, Construction of
non-isomorphic regular tournaments, Ann. Discrete Math. 52 (1992), 11-23.

\smallskip

[4] L.W. Beineke, F. Harary, The maximum number of strongly connected
subtournaments, Canad. Math. Bull. 8 (1965), 491-498.

\smallskip

[5] D.M. Berman, On the number of $5$-cycles in a tournament, in: Proc.
Sixth Southeastern Conf. Combinatorics, Graph Theory and Computing (Florida
Atlantic Univ., Boca Raton, Fla., 1975), pp. 101-108. Congress Num. XIV, Utilitas
Math., Winnipeg, Man., 1975.

\smallskip
[6] D. Burke, B. Lidick\'y, F. Pfender, M. Phillips,
Inducibility of 4-vertex tournaments,
Preprint arXiv:2103.07047v1 (2021).

\smallskip

[7] A. Brauer and I.C. Gentry, On the characteristic  roots of tournament
matrices, Bull. Amer. Math. Soc. 74 (1968), 1133-1135.

\smallskip

[8] A. Brauer and I.C. Gentry, Some remarks on tournament
matrices, Linear Algebra Appl. 5 (1972), 311-318.

\smallskip

[9] R.A. Brualdi, J. Shen, Landau's inequalities for tournament scores and
a short proof of a theorem on transitive sub-tournaments, J. Graph Theory
38 (2001), 244-254.

\smallskip

[10] U. Colombo, Sui circuiti nei grafi completi, Boll. Unione Mat. Ital.
19 (1964), 153-170.

\smallskip
[11] L.N. Coregliano and A. Razborov, On the density of transitive tournaments,
J. Graph Theory 85 (2017), 12-21.

\smallskip
[12] L.N. Coregliano, Quasi-carousel tournaments,
J. Graph Theory 88 (2018), 192-210.

\smallskip
[13] L.N. Coregliano, R.F. Parente, C.M. Sato, On the maximum density of
fixed strongly connected tournaments,
Electronic J. Comb. 26 (1), $\#$ P.1.44 (2019).

\smallskip
[14] A. Grzesik, D. Kr\'al', L.M. Lov\'asz, J. Volec,
Cycles of a given length in tournaments,
J. Combin. Theory Ser. B 158 (2023), 117-145.

\smallskip

[15]  M.G. Kendall and B. Babington Smith, On the method of paired comparisons,
Biometrika 33 (1940), 239-251.

\smallskip

[16] N. Komarov, J. Mackey, On the number of 5-cycles in a tournament,
J. Graph Theory 86 (2017), 341-356.

\smallskip

[17] G. Korvin, Some combinatorial problems on complete directed graphs,
in: Theory of Graphs (P. Rosenstiehl, Ed.), Gordon and Breach,
New York, and Dunod, Paris, 1967, pp. 197-203.

\smallskip

[18] A. Kotzig, Sur le nombre des $4$-cycles dans un tournoi,
Mat. Casopis Sloven. Akad. Vied. 18 (1968), 247-254.

\smallskip

[19] E.S. Mahmoodian, M. Mirzakhani,
Decomposition of complete tripartite graphs into $5$-cycles,
in: Combinatorics advances (Tehran, 1994). Mathematics and its
applications (C.J. Colbourn, E.S. Mahmoodian, Eds.), Springer,
Boston, MA, vol. 329 (1995), pp. 235-241; Kluwer Acad. Publ.,
Dordrecht, 1995, pp. 75-81.

\smallskip

[20] J.W. Moon, On subtournaments of a tournament,
Canad. Math. Bull. 9 (1966), 297-301.

\smallskip

[21] J.W. Moon, {\sl Topics on tournaments}, Holt, Rinehart and Winston,
New York, 1968.

\smallskip

[22] A. Moukouelle, Construction of a new class of near-homogeneous tournaments,
C.R. Acad. Sci. Paris Ser. I  327 (1998), 913-916.

\smallskip

[23] A. Moukouelle, Morphology of tournaments and distribution of 3-cycles,
PhD. Thesis, University Aix-Marseille III (1998).

\smallskip

[24] J. Plesnik, On homogeneous tournaments, Acta Fac. Rerum Natur. Univ. Comenian
Math. Publ. 21 (1968), 25-34.

\smallskip
[25] K.B. Reid and E. Brown, Doubly regular tournaments are equivalent
to skew-Hadamard matrices, J. Combin. Theory Ser. A 12 (1972), 332-338.

\smallskip

[26] K.B. Reid and L.W. Beineke, Tournaments, in: Selected topics in Graph
Theory, Vol. 2 (L.W. Beineke and R.J. Wilson, Eds.), Academic Press, New
York, 1979, pp. 169-204.

\smallskip

[27] P. Rowlison, On $4$-cycles and $5$-cycles in regular tournaments, Bull. London
Math. Soc. 18 (1986), 135-139.

\smallskip
[28] S.V. Savchenko, On $5$-cycles and $6$-cycles in regular $n$-tournaments,
J. Graph Theory 83 (2016), 44-77.

\smallskip

[29] S.V. Savchenko, On the number of $7$-cycles
in regular $n$-tournaments,
Discrete Math. 340 (2017), 264-285.

\smallskip

[30] S.V. Savchenko, On the number of $8$-cycles for two particular
regular tournaments of order $n$ with diametrically opposite local properties,
arXiv:2403.07629 (2024).

\smallskip

[31] S.V. Savchenko, On the number of $9$-cycles in a doubly-regular tournament,
(in preparation).

\smallskip

[32] S.V. Savchenko, Bernoulli numbers and $m$-cycles in regular $n$-tournaments,
Proc. Amer. Math. Soc. (will be submitted for publication).

\smallskip

[33] C. Tabib, Caract\'erisation des tournois presqu'homog\`enes, Ann. Discrete
Math. 8 (1980), 77-82.

\smallskip

[34] C. Thomassen, Hamiltonian-connected tournaments, J. Combin. Theory Series B
28 (1980), 142-163.

\end{document}